\documentclass[a4paper, 11pt]{article}  

\usepackage{graphicx} 
\usepackage{amsmath} 
\usepackage{amssymb}  
\usepackage{color}
\usepackage{hyperref}

\DeclareMathOperator*{\argmax}{arg\,max}

\title{\LARGE \bf
Dissecting Demand Response Mechanisms:\\
the Role of Consumption Forecasts and Personalized Offers\footnote{
Published in the proceedings of the 2016 Americal Control Conference (ACC 2016), Boston, USA, July 6-8 2016}
}

\author{Alberto Benegiamo\\
EURECOM and Inria, France\\
{\tt alberto.benegiamo@eurecom.fr}
\and
 Patrick Loiseau\\
 EURECOM, France\\
 {\tt patrick.loiseau@eurecom.fr}
 \and
Giovanni Neglia\\
Inria, France\\
{\tt giovanni.neglia@inria.fr}.}

\date{July 2016}

\begin{document}

\maketitle
\thispagestyle{empty}
\pagestyle{empty}

\begin{abstract}
Demand-Response (DR) programs, whereby users of an electricity network are encouraged by economic incentives to re-arrange their consumption in order to reduce production costs, are envisioned to be a key feature of the smart grid paradigm. Several recent works proposed DR mechanisms and used analytical models to derive optimal incentives. Most of these works, however, rely on a macroscopic description of the population that does not model individual choices of users. 

In this paper, we conduct a detailed analysis of those models and we argue that the macroscopic descriptions hide important assumptions that can jeopardize the mechanisms' implementation (such as the ability to make personalized offers and to perfectly estimate the demand that is moved from a timeslot to another). Then, we start from a microscopic description that explicitly models each user's decision. We introduce four DR mechanisms with various assumptions on the provider's capabilities. Contrarily to previous studies, we find that the optimization problems that result from our mechanisms are complex and can be solved numerically only through a heuristic. We present numerical simulations that compare the different mechanisms and their sensitivity to forecast errors. At a high level, our results show that the performance of DR mechanisms under reasonable assumptions on the provider's capabilities are significantly lower than those suggested by previous studies, but that the gap reduces when the population's flexibility increases. 
\end{abstract}

\section{Introduction}
Demand Response (DR hereinafter) programs are envisioned to be a key 
feature of the Smart Grid paradigm~\cite{Alb08}. By means of economic incentives (discounts or penalties), DR schemes encourage users to rearrange their consumption in response to the network state, thus mitigating the grid overload and driving wholesale prices down. 

Several analytical models are available in the literature, which describe and quantify the effects of DR mechanisms. Whatever their specifics are, these schemes need to model how users react to the incentives. Ideally the models should capture the most realistic features of a practical DR mechanism while maintaining tractability.

Among these contributions, the authors of \cite{Mung Chiang} study how an energy provider should select time-dependent discounts to minimize its production costs. They assume that the percentage of users who shift their consumption from slot $i$ to slot $j$ is a decreasing function of the temporal distance between slots $i$ and $j$ and a concave and increasing function of the discount offered in slot $j$ ($R_j$), \emph{independent} from discounts in other slots. 
The same user's model as in \cite{Mung Chiang} is adopted also in \cite{Pan14}, where the optimization problem is extended in order to account for battery storages and distributed renewable sources available into a specific microgrid.
Authors of \cite{Subra13} propose a day ahead pricing scheme which maximizes the provider's profitability and capacity utilization. Users are assumed to reschedule their consumption by comparing the utility $v_i$ they get by scheduling a task in each timeslot $i$; therefore they allocate their consumption proportionally to these utilities, i.e., they consume a fraction $\frac{v_i}{\sum_{j=1}^T v_j}$ of their total energy demand in timeslot $i$. The resulting optimization problem is non convex but some relaxation techniques are introduced, which allow one to calculate a solution within a reasonable amount of time.
In \cite{Yang13}, a more realistic model is proposed where each user  first calculates the welfare (defined as utility minus time-dependent cost) she gets from consuming electricity in each of the possible timeslots, and then allocates all the consumption to the slot returning the largest welfare. 
As we show below (see Sec.~\ref{sec.broadcast}) this model can lead to a much more complex optimization problem than the one presented in \cite{Yang13}. 
Finally, the authors of \cite{Song14} propose a full-fledged game theoretical model, but their results hold only if users experience a large number of interactions without any change in the system.

We claim that these studies rely on too strong assumptions, which jeopardize their usability for practical purposes. Interestingly, we observe that the assumptions are sometimes hidden in the macroscopic models the papers start from. In particular in this paper we focus on \cite{Mung Chiang} and show that its model requires personalized offers and a very precise forecast of the baseline consumption of each user. 
The implementation of these features may require potentially significant costs in terms of communication, measurement and computation infrastructure.
Besides highlighting these implicit requirements in the analytical framework in \cite{Mung Chiang} (and then also in \cite{Pan14}), we explore their potentials considering four DR mechanisms with different levels of complexity:
\begin{enumerate}
 \item the \emph{base} mechanism corresponds to an optimization problem similar to the one considered in~\cite{Mung Chiang}, it requires personalized offers and individual consumption forecasts; the energy production cost is optimized over the discount values, each of which is offered to a given fraction of the population,
 \item the \emph{optimized} mechanism takes full advantage of personalized offers and consumption forecasts by minimizing the cost over both the discount values and the population fractions to which the discounts are offered,
 \item the \emph{robust} mechanism relies on personalized offers, but does not need individual consumption forecasts, 
 \item finally the \emph{broadcast} mechanism (analogous to that in~\cite{Yang13}) needs neither of the two features.
 \end{enumerate}
Interestingly, contrarily to prior studies, we find that the cost-minimization problems resulting from our DR mechanisms are not convex (even for the base mechanism). Nevertheless, simple heuristics can identify (potential) minima in a reasonable amount of time in realistic scenarios. 
Then, our numerical results show that the simpler robust and broadcast mechanisms achieve significantly lower cost reductions than the optimized mechanism, which is difficult to implement, but that the gap reduces when the population's flexibility increases. 

The paper is organized as follows. In Sec.~\ref{sec.macroscopic} we discuss how the macroscopic models considered in~\cite{Mung Chiang, Pan14, Subra13} hide some implicit assumptions about the user rationality or about the interactions between the provider and the user. We define our microscopic model in Sec.~\ref{s:microscopic} and then describe different DR mechanisms and their corresponding optimization problems in~Sec.~\ref{sec.mechanisms}. Finally, we evaluate their performance numerically in a realistic scenario in Sec.~\ref{sec.results}. 

Due to space constraints proofs, examples and some additional numerical results are in the companion technical report~\cite{Ben15tr}.

\section{Pitfalls when Starting from Macroscopic Models}\label{sec.macroscopic}
In this section, we describe in more detail the macroscopic models proposed in the literature for day-ahead price optimization. Consider a finite time horizon discretized in a set $\mathcal T$ of $N$ timeslots and a large population $\mathcal{S}$ of users. The baseline aggregate energy consumption in slot $j$ is denoted by $E^0_j$. 

The energy provider charges a flat rate $B$, but it can offer discount rates to incentivize the users to move some of their consumption so as to reduce the energy production cost. Due to consumption shifts, the actual aggregate consumption in time slot $j$ is $E^1_j$. 
Observe that a usual assumption in the literature (including the papers mentioned above) is that the introduction of a DR scheme neither reduces nor increases users' demand; it merely rearranges users' consumption in a more cost effective way, so that
\begin{equation}\label{conservation}
\sum_{j=1}^N E_j^0 = \sum_{j=1}^N E_j^1.
\end{equation}
We denote  the amount of consumption shifted from slot $j$ to slot $i\neq j$ as $E_{j \rightarrow i}$, and the amount of consumption the users refuse to shift away from $j$ as $E_{j \rightarrow j}$. Then we have
\[E^1_i= E_i^0 + \sum_{z=1}^N E_{z \rightarrow i} - \sum_{k=1}^N E_{i \rightarrow k}.\]

We now start to further detail the model considering some specific assumptions made in previous works.
In~\cite{Mung Chiang} and~\cite{Pan14}, the electricity provider offers an energy price discount $R_i \ge 0$ in each slot $i$. The users are assumed to react to these incentives by shifting a fraction of their baseline consumption from slot $j$ to slot $i$ ($|j-i|$ slots away) according to the following formula: 
\begin{equation}
\label{e:sensitivity}
E_{j \rightarrow i} = E^0_j S_j(R_i,|j-i|).
\end{equation} 
$S_j(R_i,|j-i|)$ is called the aggregate sensitivity function and is increasing in the discount $R_i$ and decreasing in the temporal shift $|j-i|$, in order to take into account the user discomfort.

The  provider selects the vector of discounts $\mathbf R$ in order to minimize its total cost, equal to the sum of the electricity generation costs and the loss of revenues due to the discounts. In particular the optimization problem considered in~\cite{Mung Chiang} is the following: 
\begin{eqnarray}
\min_{\mathbf{R}} &  \sum_i \sum_{j \neq i} R_i E_{j \rightarrow i} +\sum_i c_i \big(E_i^1  \big) \label{Mech_1:1}\\
\text{s.t.} & 0 \leq R_i \leq B \quad \forall i=1,\ldots N, \label{Mech_1:2}
\end{eqnarray}
where $c_i(\cdot)$ is the cost of electricity production in slot $i$. 
Eq.~\eqref{Mech_1:2} guarantees that discounts $\mathbf{R}$ are non negative and smaller than the flat rate $B$, so that the money stream goes toward the provider.

As it often happens, the devil is hidden in the details, and in this case in Eqs.~\eqref{e:sensitivity} and~\eqref{Mech_1:1}. 
Our first remark is that the cost of lost revenues $ \sum_i \sum_{j \neq i} R_i E_{j \rightarrow i}$ in Eq.~\eqref{Mech_1:1} implicitly assumes the possibility to reward only the consumption actually shifted from $j$ to $i$, i.e., $E_{j \rightarrow i}$, but this quantity cannot be directly measured. The actual consumption $E^1_i$ can be measured, and then $E_{j \rightarrow i}$ can be quantified provided that we have good estimates of the sensitivity function $S_j(R_i,|j-i|)$ and of the baseline consumption $E^0_i$. Let us assume for a moment that $S_j(R_i,|j-i|)$ is known from historical data and that the aggregate baseline consumption may be predicted with a reasonably high level of accuracy on a large set of users. Then it seems possible to solve the macroscopic problem in Eqs.~\eqref{Mech_1:1} and~\eqref{Mech_1:2}, but we need to consider also what should happen at the microscopic scale. While the estimates for the aggregate baseline consumption can be adequately precise, finally the billing is done at the user's granularity and each user expects to receive the price discount corresponding to the energy consumption she actually moved. 
If the energy bill's reduction does not correspond to her forecast, the user is likely to opt out of the program (in particular if she has experienced underpayments) or to reduce her efforts and milk occasional discounts.
It appears then that Eq.~\eqref{Mech_1:1} implicitly requires very precise predictions of \emph{individual} consumptions.

We now observe that the form of the sensitivity function $S_j(R_i,|j-i|)$ in Eq.~\eqref{e:sensitivity} indicates that the amount of energy shifted from $j$ to $i$  depends on the discount $R_i$ but not on the other discounts.
We can then ask ourselves which individual decisions may lead to this aggregate behavior, an issue ignored both in~\cite{Mung Chiang} and~\cite{Pan14}. As long as a rational individual is offered two different discounts $R_i$ and $R_k$, it seems natural that her decision to move some consumption from $j$ to $i$ or from $j$ to $k$ or to keep it in $j$ will take into account both the discounts. 
To stress the point, consider a case when both $S_j(R_i,|j-i|)$ and $S_j(R_k,|j-k|)$ are positive, but moving the consumption from $i$ to $k$ is both less inconvenient (i.e.,~$|j-k|<|j-i|$) and more rewarding (i.e.,~$R_k>R_i$).  There is then no reason why the user would move consumption to $i$. The conclusion is the same for all the users 
and then we should have $E_{j \rightarrow i}=0$ at the aggregate level, in contradiction with Eq.~\eqref{e:sensitivity}.
We can then conclude that the expression of the sensitivity function in Eq.~\eqref{e:sensitivity} is not suited to model the situation when a user is offered two or more rewards, but it can capture the case when the user decides between moving from $j$ to $i$ in exchange of a discount $R_i$ or staying in $j$. 
If every user is offered a single discount to move to a given slot, but different users can receive different offers, then Eq.~\eqref{e:sensitivity} can reasonably describe the macroscopic effect of such personalized offers.
The details are described in Sec.~\ref{s:base}, here we only highlight that Eq.~\eqref{e:sensitivity} requires then that the electricity provider i) calculates an offer for each user, ii) communicates individually to the user, iii) considers the individual offer when billing the user.
This is clearly more demanding than simply advertising to the whole population the same set of discounts.

We observe that the equivalent sensitivity function considered in~\cite{Subra13} poses similar problems. Using our notation, we have $E_{j \rightarrow i}=\frac{v_i}{\sum_{k\in \mathcal T} v_k} E^0_j$, where $v_i$ is the net utility a user gets by consuming electricity in slot $i$ and can be a function of the timeslot itself and of the discount $R_i$. This formula tries to capture the effect of the whole set of discounts, but it is not clear again what is the underlying user's model: if slot $i$ has a larger utility than slot $k$ ($v_i>v_k$), why should the user consume in $k$?

\section{Starting from a Microscopic Model}
\label{s:microscopic}
In the previous section we made the point that, while aggregate population models may be convenient, it is necessary to explicitly consider the microscopic level: how the user takes the decisions and how the provider and the user are supposed to interact. 
In this paper we follow the opposite path in comparison to the existing works mentioned: we move from the microscopic level to the macroscopic one. 
In particular, in this section, we start from a clear model of rationality for the single user and then move to describe how aggregate quantities can be derived. 

Each user $u$ has a baseline energy consumption $\{e_j^{0u}\}_{j=1,\ldots N}$
that leads to the aggregate baseline consumption $E^0_j=\sum_{u \in \mathcal S}  e_j^{0u} $ for $j=1,\ldots N$.
We assume in what follows that users are homogeneous, i.e., 
\begin{equation}
\forall \; u, \quad e_j^{0u}=e_j^0 \quad j=1,\ldots N.
\end{equation}
In~\cite{Ben15tr} we show how the DR mechanisms perform when this assumption does not hold.

User $u$ is characterized by a private type $\underline{D}^u=\{d^u_{j \rightarrow i}\}_{j,i=1, \ldots N}$ where $d^u_{j \rightarrow i}$ indicates the discomfort due to shifting one unit of consumption from timeslot $j$ to timeslot $i$. We assume that discomforts are expressed in monetary units; and that, $\forall u \in \mathcal{S}$,
\begin{equation}
d^u_{j \rightarrow j} =0 \;\; \textrm{ and } \;\; d^u_{j \rightarrow i} >0, \;\; \forall j, i \neq j,
\end{equation}
i.e., there is a strictly positive discomfort if and only if consumption is shifted from its original timeslot. The provider does not know the private type $\underline{D}^u$ of each user $u$: from its point of view, each discomfort $d^u_{j \rightarrow i}$ is drawn from a known, continuous distribution $F_{j \rightarrow i}$ on $[0, \alpha_{j,i}]$ (where possibly $\alpha_{j,i}=+\infty$). 
Discomforts of distinct users are independent but note that we do not assume that, for a given user, the discomforts $\{d^u_{j \rightarrow i}\}_{j,i=1, \ldots N}$ are mutually independent. 

\subsection{Rational Users}\label{sec.rational}

We assume that a user simply chooses the option that maximizes her utility.
In particular let $\mathcal{T}^u_j$ be a set of timeslots the user could move the baseline consumption $e^{0u}_j$ to in exchange for different discounts $\mathcal{R}_j^u=\{R^u_{j \rightarrow k}\ge0, k \in \mathcal{T}^u_j\}$. The set pair $(\mathcal{T}^u_j,\mathcal{R}_j^u)$ defines the offer user $u$ receives for timeslot $j$.
The set of options includes the possibility to keep the consumption in $j$, i.e., $j \in \mathcal{T}^u_j$.
A rational user maximizes her utility by scheduling her consumption $e_j^{0u}$ to a timeslot
\begin{equation}
\label{e:argmax}
i^* \in \argmax_{k \in \mathcal{T}^u_j}\left\{R^u_{j \rightarrow k} - d^u_{j \rightarrow k}\right\}.
\end{equation}
We assume that if two or more timeslots are equally palatable, the whole consumption is shifted to only one of them, picked at random with equal probability.

\subsection{Aggregation}
We observe that the quantities $d^u_{j \rightarrow k}$ in Eq.~\eqref{e:argmax} are random variables, then two different users could take different decisions while confronted with the same offers. 
The aggregate consumption $E^1_i$, for $i\in \mathcal T$, would then be a random variable. Here we assume (as it is implicit in the other works) that we always work with large sets of the population so that the variability can be neglected by approximating actual random quantities with their expected values. In particular, if a subset $\mathcal{Q}$ containing a fraction $q$ of the population receives an offer $(\mathcal{T}_j, \mathcal{R}_j)$, the corresponding consumption shifted from $j$ to a time slot $i$, denoted as $E_{j \rightarrow i}^{\mathcal Q}$ will be
\begin{equation}
\label{e:aggregated}
 E_{j \rightarrow i}^{\mathcal Q} = q E^0_j \textrm{ Prob}\left(i  \in \argmax_{k \in \mathcal{T}_j}\left\{R^u_{j \rightarrow k} - d^u_{j \rightarrow k}\right\}\right),
\end{equation}
if the probability that a user has two or more equally palatable timeslots is zero.
When discomforts are continuous random variables (as we consider in this paper), this is always the case if each user receives only one offer (the first three mechanisms introduced below) or if the discomforts $\{d^u_{j \rightarrow i}\}_{j,i=1, \ldots N}$ are mutually independent.
In Sec.~\ref{sec.broadcast}, we discuss how Eq.~\eqref{e:aggregated} should be modified if this probability is not zero.
We denote $\textrm{ Prob}\left(i  \in \argmax_{k \in \mathcal{T}_j}\left\{R_{j \rightarrow k} - d^u_{j \rightarrow k}\right\}\right)$ simply as $P_{j \rightarrow i} (\mathcal R_j^u, \mathcal T_j)$.

\section{DR mechanisms}\label{sec.mechanisms}
Under different assumptions on the provider's capabilities, we introduce different demand response mechanisms based on the microscopic model above, which are therefore practically implementable. We introduce and study the corresponding optimization problems. 

We start by the base mechanism that leads to the same aggregate optimization problem considered in~\cite{Mung Chiang, Pan14}.

\subsection{Base mechanism} \label{s:base}
This mechanism requires that the energy provider can manage personalized offers to its customers and moreover that it has perfect knowledge (or very precise estimates) of the baseline consumption of each user.

The population is segmented into $N^2$ disjoint subsets $\mathcal{Q}_{j \rightarrow i}$, for $j, i \in \mathcal T$, respectively including a fixed fraction $q_{j \rightarrow i}$ of the population.
Each user in  $\mathcal{Q}_{j \rightarrow i}$ is simply offered to move her baseline consumption in slot $j$ ($e^0_{j}$) to slot $i$ in exchange for a price discount $R_i$.

The total consumption that is shifted from $j$ to $i$ is then 
\[E_{j \rightarrow i} = q_{j \rightarrow i} E^0_j \textrm{ Prob}\left(R_i - d^u_{j \rightarrow i} > 0\right)\]
as it can be obtained from Eq.~\eqref{e:aggregated}, taking into account that in this case $\mathcal{T}_j=\{j,i\}$ and $R_j-d^u_{j \rightarrow j}=0$. We observe that the probability appearing on the right-hand side only depends on the reward $R_i$ and on the random variable  $d^u_{j \rightarrow i}$. If the discomfort is only a function of the temporal distance $|j-i|$, then the  sensitivity function (the ratio of people who move from $j$ to $i$) has the same properties than in~\cite{Mung Chiang}, in particular:
\[S_j(R_i,\vert i-j \vert)=q_{j \rightarrow i}P_{j \rightarrow i}(R_i),\]
where for $P_{j \rightarrow i}(\cdot)$ we have made explicit the only variable it depends from.

As we discussed in Sec.~\ref{sec.macroscopic}, because the provider knows exactly the consumption shifted from each user, it can formulate the optimization problem~(\ref{Mech_1:1}-\ref{Mech_1:2}).
In \cite{Mung Chiang} it is stated that the problem is convex if i)
the productions costs $c_j(\cdot)$ are continuous piecewise linear and increasing and ii) the discomfort distributions $F_{j \rightarrow i}(\cdot)$ are continuous and concave. We show in~\cite{Ben15tr} that this is not the case by providing a counterexample. Stronger hypotheses are required for the problem to be concave, as for example the linearity of the discomfort functions.

In particular in \cite{Mung Chiang} the numerical evaluation considers 
\[S_{j}(R_i,\vert i-j\vert)= \frac{1}{\sum_{k=1}^{N}\frac{1}{(\vert k-j \vert+1)}}\frac{R_i}{B \cdot \left(\vert i-j\vert+1\right)},\] 
that leads us to consider 
\begin{equation}\label{Mung size}
q_{j \rightarrow i}=\frac{\frac{1}{\vert i-j \vert+1}}{\sum_{k=0}^N \frac{1}{\vert k-j \vert +1}}, \;\;\; P_{j \rightarrow i}(R_i)= \frac{R_i}{B}.
\end{equation}
This particular expression for $P_{j \rightarrow i}$ can be obtained if $d_{j \rightarrow i}^u$ is a uniform random variable with support in $[0,B]$.
The numerical results for the base mechanism in Sec. \ref{sec.results} are obtained considering the same expression for the fractions $q_{j \rightarrow i}$.

Due to the non-convexity of the optimization problem~(\ref{Mech_1:1}-\ref{Mech_1:2}) we cannot use one of the classic algorithms for convex optimization. For the results shown in section \ref{sec.results} we have adopted instead a multi-start approach: we have generated random starting points uniformly distributed in the problem domain and we have run per each point a descendent algorithm which converged on a local minimum; the optimal offers are therefore those returning the smallest cost among these minimizers. This approach does not guarantee convergence to the global optimum but its reliability can be improved by increasing the number of starting points.

\subsection{Optimized Mechanism}
We have now understood which DR mechanism can lead to the optimization problem~\eqref{Mech_1:1}, but now that we look at its implementation at the microscopic level 
and the need for personalized offers,
some specifics of the base mechanism look arbitrary and unjustified.
For example, given that discounts are not broadcast but each user receives an individual offer, 
why should the discounts offered to the two disjoint sets of users $\mathcal Q_{j \rightarrow i}$ and $\mathcal Q_{k \rightarrow i}$ be equal to the same value $R_i$? 
It is clear that the energy provider can further reduce the cost if it can independently choose $R_{j \rightarrow i}$ and $R_{k \rightarrow i}$.
Moreover, there is no reason to think that the size of the sets $\{\mathcal Q_{j \rightarrow i}\}$ should be fixed, the fractions $\{q_{j \rightarrow i}\}$ can also be optimization variables.

We allow the provider to take advantage of these additional degrees of freedom that---we repeat---do not impose any additional requirement to the system. We call this new DR mechanism optimized.
 The load $E_{j \rightarrow i}(R_{j \rightarrow i},q_{j \rightarrow i})$ rescheduled from $j$ to $i$ is now
$E_{j \rightarrow i}(R_{j \rightarrow i},q_{j \rightarrow i})= q_{j \rightarrow i} P_{j \rightarrow i}(R_{j \rightarrow i}) E_j^0$
and the cost minimization problem becomes:
\begin{eqnarray}
\min_{\mathbf{R},\mathbf{q}} & \!\!\!\!\! \text{  cost}_{\text{opt.}}(\mathbf{R},\mathbf{q})= \sum_{i} \!\sum_{z\neq i}\!\! R_{ z \rightarrow i} E_{z \rightarrow i} \! +\! \sum_i c_i \!\big(\!E_i^1 \! \big) \label{Mech_2:1}\\
\text{s.t.} & 0 \leq R_{ z \rightarrow i} \leq B \quad \forall z,i=1,\ldots N \label{Mech_2:2}\\
 &0 \leq q_{z \rightarrow i} \leq 1, \quad z,i=1,\ldots N \label{Mech_2:3}\\
 &\sum_{i} q_{z \rightarrow i} \leq 1, \quad \forall z=1 \ldots N. \label{Mech_2:4}
\end{eqnarray}
Eq. \eqref{Mech_2:2} guarantees that discounts $\mathbf{R}$ are non negative and smaller than the flat rate $B$, Eq.~\eqref{Mech_2:4} is a consequence of the fact that each user receives at most one offer for its baseline consumption in a given slot.

The optimization problem (\ref{Mech_2:1}-\ref{Mech_2:4}) can be solved with the same heuristic proposed for problem (\ref{Mech_1:1}-\ref{Mech_1:2}).

\subsection{Robust Mechanism}
The optimization problems~(\ref{Mech_1:1}-\ref{Mech_1:2}) and (\ref{Mech_2:1}-\ref{Mech_2:4}) assume that the provider has perfect knowledge of each user's baseline consumption, so that it can correctly identify the consumption shifted and reduce accordingly the energy bill.
This assumption is probably unrealistic. 
If the provider does not have such capability, then it can offer the user a discount for all the consumption in a given timeslot $i$ and not just for the consumption moved to $i$. 
The population is then divided into $N$ subsets $\mathcal{Q}_i$, each containing a fraction $q_i$ of the users. All users in $\mathcal{Q}_i$ receive one and only one offer: they are encouraged to shift their consumption from any timeslot in the time horizon to timeslot $i$ and they get the discount $R_i$ for all the electricity consumed in $i$, including the one originally in $i$. 

We call this scheme robust, because it does not rely on estimates of individual consumption.
It is clearly simpler than the previous two, because the provider needs only to measure the amount of consumption in $i$ for the users who got the offer and to bill them accordingly.

The load $E_{j \rightarrow i}(R_{i},q_{i})$ shifted from $j$ to $i$ is
$E_{j \rightarrow i}(R_{i},q_{i})= q_{i} P_{j \rightarrow i}(R_{i}) E_j^0$.
Note that users in $\mathcal{Q}_i$ have no interest to move their baseline consumption away from $i$, then  $E_{i \rightarrow i} = q_i E^0_i$. 
The robust mechanisms leads to the following optimization problem:
\begin{eqnarray}
\min_{\mathbf{R},\mathbf{q}} &\text{  cost}_{\text{rob.}}(\mathbf{R},\mathbf{q})=\sum_{i,z} R_i E_{z \rightarrow i} +\sum_i c_i \big(E_i^1 \big) \label{Mech_3:1}\\
\text{s.t.} & 0 \leq R_{i} \leq B \quad i=1,\ldots, N \label{Mech_3:2}\\
& 0 \leq q_{i} \leq 1, \quad i=1,\ldots, N \label{Mech_3:3}\\
& \sum_{i=1}^N q_{i} \leq 1.\label{Mech_3:4} 
\end{eqnarray}
Note that in Eq.~\eqref{Mech_3:1} the first sum includes also $E_{i \rightarrow i}$ because all the final consumption in $i$ from the users in $\mathcal Q_i$ is paid at a discounted price. The term does not appear in Eq.~\eqref{Mech_1:1} and Eq.~\eqref{Mech_2:1}.

The optimization problem (\ref{Mech_3:1}-\ref{Mech_3:4}) can be solved with the same heuristic proposed for problem (\ref{Mech_1:1}-\ref{Mech_1:2}).

\subsection{Broadcast Mechanism}\label{sec.broadcast}

In the three mechanisms introduced above, the provider makes personalized offers to users in selected fractions of the population. This may not always be possible (due to the complexity it introduces for instance in billing) or desirable (for perceived fairness issues). Our last mechanism, which is the simplest (in its definition), does not assume personalized offers. The provider selects a single vector $\mathbf{R}$ of discounts for every time slot and broadcasts these discounts to all users (hence the name broadcast mechanism). Users then re-arrange their demand and pay the discounted price for their demand in each slot (hence this mechanism also does not rely on the need to estimate shifted demand). 

As explained in Sec.~\ref{sec.rational}, each individual user moves her demand from slot $j$ to a slot (potentially $j$ itself) that maximizes her net utility (discount minus discomfort). Recall that if several slots give equal net utility, the user chooses one of them randomly. 

Until now, we have not made any assumption on the possible correlation of a given user's discomforts. This is because, in the previous three mechanisms, each user was receiving only one offer.
In the broadcast mechanism, each user has several offers to compare to decide on his new demand schedule, we therefore need to describe the discomfort correlations. 

Let us consider now the particular case when two slots, say $h$ and $k$ may appear equally attractive to a user, i.e.~$R_h-d_{j \rightarrow h}=R_k-d_{j \rightarrow k}$.
If we assumed that, for each user, the discomforts $\{d_{j \rightarrow i}\}_{\{i, j =1, \cdots, N\}}$ were mutually independent, this event would have probability zero according to our assumption on $F_{i\rightarrow j}$, and therefore it would not appear at the aggregate level. As a result, the aggregate demand moved from $j$ to $i$ would be
$E_{j \rightarrow i}(\mathbf{R}) = P_{j \rightarrow i}(\mathbf{R}) E_j^0$,
where
\begin{equation}
\label{eq.Pbroadcast} P_{j \rightarrow i}(\mathbf{R}) = \textrm{Pr}\left(R_i -d_{j \rightarrow i} \geq \max_{k \neq i} \left\{R_k- d_{j \rightarrow k} \right\} \right). 
\end{equation}

However, rather than making the above independence assumption, we prefer to assume that the discomforts have the form $d_{j \rightarrow i} = \beta_j \vert i-j \vert^{t_j}$, where $t_j$ is a constant independent of the user and $\beta_j$ is a random variable with concave Cumulative Distribution Function (CDF) $F_j(\cdot)$. 
This model describes a symmetric delay sensitivity of users (users are indifferent between moving two hours earlier or two hours later) while keeping the flexibility of users having a different flexibility of demand of different times (since $\beta$ and $t$ are indexed by the origin timeslot $j$); but it also introduces correlations between the discomforts of a user. As a result, the fraction of demand shifting from  $j$ to $i$ is 
\begin{equation}
\label{eq.Pbroadcast2} P_{j \rightarrow i}(\mathbf{R}) \! =\! \frac{ \textrm{Pr}\left(R_i -d_{j \rightarrow i} \geq \max_{k \neq i} \left\{R_k- d_{j \rightarrow k} \right\} \right)}{1+\mathbf{1}_{R_{2j-i}=R_i}},
\end{equation}
rather than \eqref{eq.Pbroadcast}. The denominator in \eqref{eq.Pbroadcast2} accounts for cases when a slot other than $i$ (which has to be $2j-i$) gives equal net utility \emph{for all users}. 
The broadcast mechanism then leads to the following optimization problem: 
\begin{eqnarray}
\min_{\mathbf{R}} & \text{  cost}_{\text{brd.}}(\mathbf{R})=\sum_z \sum_i R_i E_{z \rightarrow i} +\sum_i c_i \big(E_i^1 \big) \label{Mech_4:1}\\
\text{s.t.} & 0 \leq R_i \leq B \quad \forall i=1,\ldots N. \label{Mech_4:2}
\end{eqnarray}

Unfortunately, due to  indicator function in Eq.~\eqref{eq.Pbroadcast2}, the cost function~\eqref{Mech_4:1} of the broadcast mechanism is not continuous, even in very simple scenarios with continuous production costs~\cite{Ben15tr}. Discontinuity arises also in the macroscopic model in~\cite{Yang13}, but it seems to have been ignored.

In practice, we solve problem~(\ref{Mech_4:1}-\ref{Mech_4:2}) using the same heuristic proposed for the previous problems, but we work on a continuous and smooth approximation of the cost function.

\subsection{Ranking DR mechanisms}
\label{s:ranking}
In~\cite{Ben15tr} we prove that the optimized mechanism always performs better than the base and robust mechanisms, but the ranking cannot be extended further.

\section{Numerical Results}\label{sec.results}
In this section we evaluate the performance of the different DR mechanisms in the realistic scenario considered in  \cite{Mung Chiang} and based on energy data about the Ontario province in Canada. In particular, the baseline consumption $\mathbf{E_0}$ and the cost functions  are estimated from the IESO energy portfolio~\cite{IESO}, consisting of nuclear plans, hydro gas powered stations and renewable and from typical costs associated to these energy sources.
We considered  the flat rate $B = \text{110\$/MWh}$ and the timeslot-independent piecewise-linear cost function $c(\cdot)$ with derivative: $c'(E)= \$ 10$ for $E \leq C_1$, $c'(E)=\$ 72.46$  for $E \in (C_1\,,\,C_2)$ and $c'(E)= \$ 91$ for $E \geq C_2$,
where $C_1 = 16.3 \text{ GWh}$ and $C_2 = 17.9 \text{ GWh}$ represent respectively the base to intermediate load capacity and intermediate to peak load capacity. 
We assume that discomforts take the form $d_{j \rightarrow i}= \beta_j \vert i-j \vert$, where $\beta_j$ is an exponential random variable with cdf $F_j(R;\beta)=1-e^{\frac{\beta}{\mu}}$. $\mu$ is a parameter representing the population's flexibility. The larger it is, the smaller (in a stochastic order sense) the discomfort of the users to shift their consumption.

\begin{figure}
\centering
\includegraphics[scale=0.8]{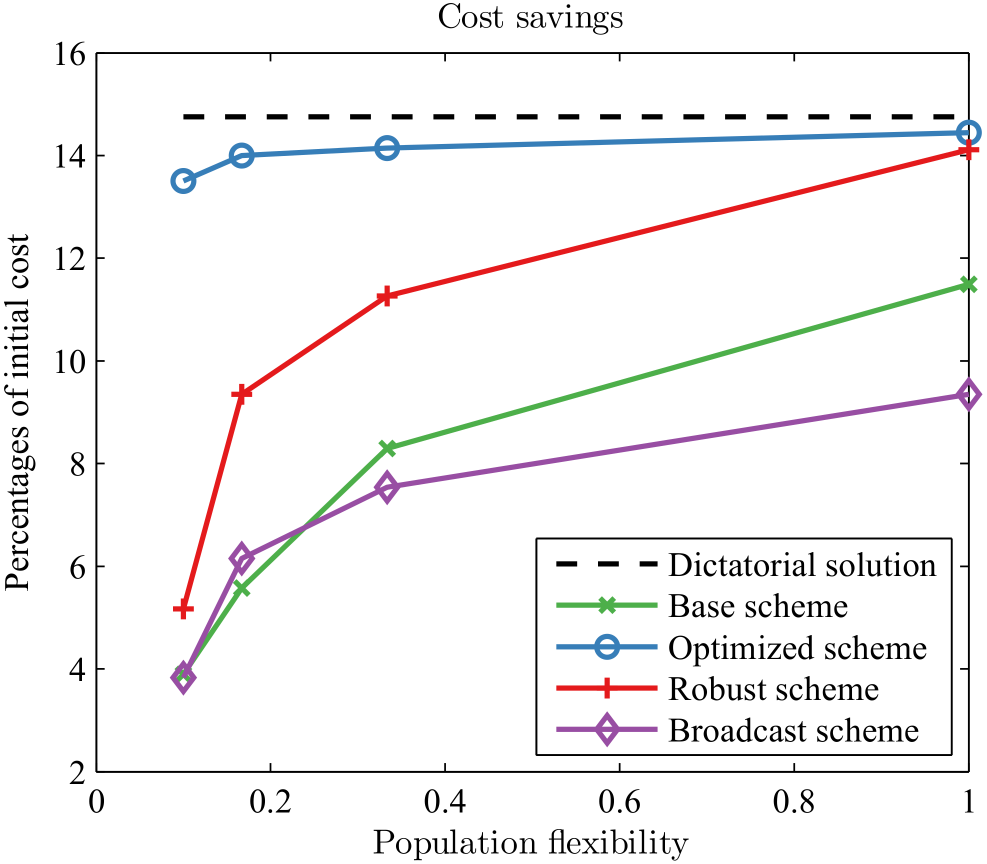}
\caption{Cost savings normalized to the initial cost, for various flexibility parameters $\mu$.}
\label{savings}
\end{figure}

In Fig. \ref{savings}, we report the cost savings of the DR schemes, normalized to the initial cost, for four different values of the $\mu$ parameter: $\frac{1}{10}$, $\frac{1}{6}$, $\frac{1}{3}$, $1$. The dashed line represents the saving which could be achieved if users' demand could be rearranged at the provider's will without providing any discount (we indicate it as the dictatorial solution).
Consistently with the results in Sec.~\ref{s:ranking}, the optimized mechanism returns larger savings than the robust and the base ones. Interestingly, the robust mechanism performs consistently better than the base one despite the fact that it does not require the ability to estimate the demand shifted and it therefore ``wastes'' some discount by giving it to demand that was already scheduled in a given timeslot in the baseline demand. 
Moreover, as the population flexibility increases, the savings gap between the optimized scheme and the robust mechanism reduces, the latter being effectively close to exploiting all the population's flexibility.

\begin{figure}
\centering
\includegraphics[scale=0.8]{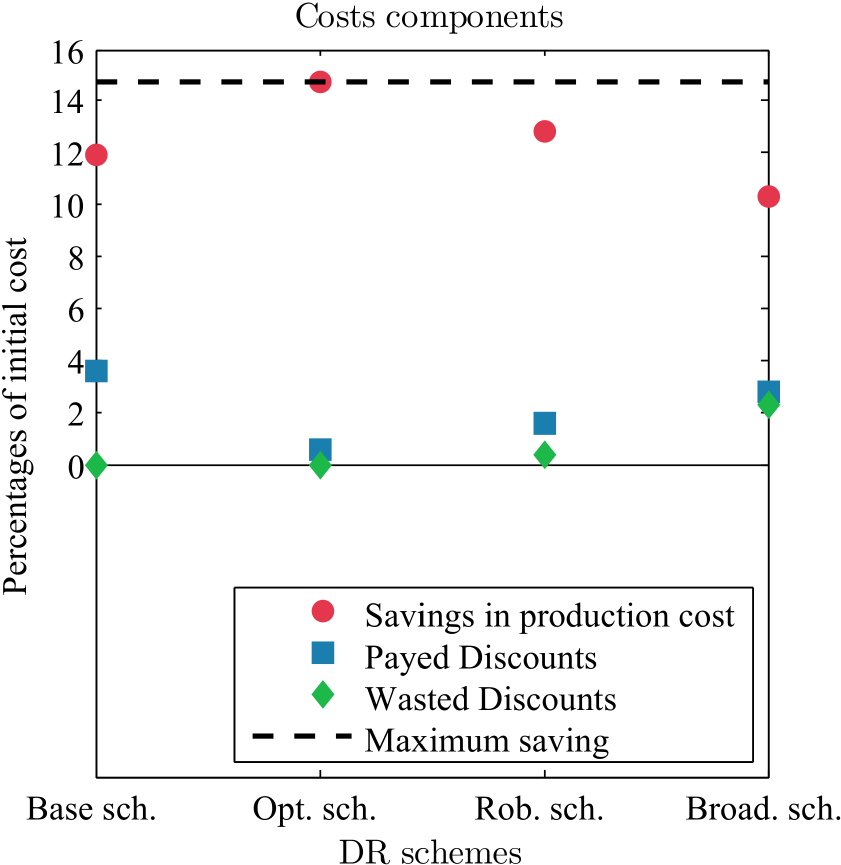}
\caption{Analysis of the components of the cost savings. All the quantities are normalized to the initial cost.}
\label{Cost analysed}
\end{figure}

In Fig.~\ref{Cost analysed}, we focus on the case $\mu = \frac{1}{3}$ and analyze the components of the cost for each DR mechanism. Fig.~\ref{Cost analysed} confirms that the optimized scheme provides the largest savings as it can minimize the production cost while paying the smallest amount of discounts. We indicate with wasted discounts the amount of discounts paid to consumption that would in any case have been scheduled in that timeslot. The base and optimized mechanisms do not waste any discount, while the robust mechanism and the broadcast scheme do, as they provide the discount $R_i$ to all the electricity consumed in $i$, including the part of $E_i^0$ that remains in $i$.

\section{Conclusions}
In this paper, we have shown that macroscopic descriptions of DR mechanisms can hide important assumptions that can jeopardize the mechanisms' implementation. For this reason,  our proposal moved from a microscopic description that explicitly models each user's decision. We have then introduced four DR mechanisms with various assumptions on the provider's capabilities. 
Interestingly, contrarily to previous studies, we find that the optimization problems that result from our mechanisms are complex and can be solved numerically only through a heuristic. 
Moreover, our results show that the performance of DR mechanisms under reasonable assumptions on the provider's capabilities are significantly lower than those suggested by previous studies, but that the gap reduces when the population's flexibility increases.

\section{Acknowledgements}
This work was partly funded by the French Government (National Research Agency, ANR) through the ``Investments for the Future'' Program reference \#ANR-11-LABX-0031-01.

\end{document}